\documentclass[11pt]{article}

\usepackage[T1]{fontenc}
\usepackage{microtype}
\usepackage{amsmath,amssymb,amsthm}
\usepackage{mathtools} % for math symbols
\usepackage{amsmath}
\usepackage{amssymb} % for \End
\usepackage[margin=1in]{geometry}
\usepackage{enumitem}
\usepackage{hyperref}

\numberwithin{equation}{section}

\theoremstyle{plain}
\newtheorem{theorem}{Theorem}[section]
\newtheorem{proposition}[theorem]{Proposition}
\newtheorem{lemma}[theorem]{Lemma}

\theoremstyle{definition}
\newtheorem{definition}[theorem]{Definition}

\theoremstyle{remark}
\newtheorem{remark}[theorem]{Remark}

\usepackage{bookmark} % safer bookmarks (recommended)
\begin{document}

\title{Deformations of Jordan Algebras via the Jordan Defect\\
An Explicit Low--Degree Deformation Complex}

\author{
Vincent E.~Coll, Jr.\\
Department of Mathematics, Lehigh University\\
Bethlehem, PA 18015\\
\texttt{vec208@lehigh.edu}
}

\date{}
\maketitle

\begin{abstract}
Over a field of characteristic $0$ we give a concrete, computation--ready description
of Jordan algebra structures and their low--order deformation theory.
The Jordan identity is quartic in the elements and cubic in the multiplication,
and in characteristic $0$ it is equivalent to its standard four--variable polarization.
We encode this polarization as a cubic map in the product~$\mu$, called the
\emph{Jordan defect} $J(\mu)$.
Linearizing this defect yields an explicit low--degree deformation complex
\[
C^1(J)\xrightarrow{\;\delta_\mu\;} C^2(J)\xrightarrow{\;d_\mu\;} C^3(J),
\]
whose second cohomology classifies infinitesimal deformations modulo equivalence
and whose obstruction space
\[
\mathrm{Obs}^3_\mu := C^3(J)/\operatorname{im}(d_\mu)
\]
contains the primary obstruction to extending such deformations.
We emphasize that this construction captures only the low--degree part of the
operadic deformation theory and does not claim to produce the full governing
$L_\infty$ structure.
\end{abstract}

\noindent\textbf{Keywords:}
Jordan algebras, deformation theory, Jordan identity, cohomology, algebraic deformations.

\noindent\textbf{Mathematics Subject Classification (2020):}
17C05 (Primary), 16E40 (Secondary).

\medskip

\begin{remark}[Relation to the full operadic picture]
The deformation theory of Jordan algebras is governed in general by an
$L_\infty$ algebra arising from a cofibrant resolution of the Jordan operad.
Because the Jordan operad is not quadratic, this $L_\infty$ structure necessarily
involves higher brackets $\ell_n$ for $n\ge3$.
In this paper we do not attempt to construct the full operadic deformation complex.
Instead, we extract explicitly the first nontrivial bracket corresponding to the
cubic Jordan identity and its linearization.
This suffices to describe infinitesimal deformations and the primary obstruction,
but not higher--order obstructions.
\end{remark}

\begin{remark}[What this viewpoint adds]
Classical deformation theory of Jordan algebras guarantees, in principle, the
existence of a controlling homotopical structure and associated cohomology groups.
However, this general theory does not typically provide an explicit, computation--ready
description of the low--degree deformation equations and equivalence relations for a
given Jordan algebra.
The point of the present construction is to make this low--degree structure completely
explicit: infinitesimal deformations are characterized by a finite linear system obtained
by linearizing a single polynomial identity, and equivalence is given by an explicit
change--of--coordinates action.
As a result, deformation and rigidity questions for specific Jordan algebras reduce
directly to finite linear algebra, without invoking resolutions or higher homotopical data.
\end{remark}

% ==================================================
\section{Jordan algebras and the Jordan defect}

Let $J$ be a vector space over a field $k$ of characteristic $0$.
A commutative bilinear map
\[
\mu : \mathrm{Sym}^2 J \to J
\]
is a \emph{Jordan product} if
\begin{equation}\label{eq:jordan}
(x^2\circ y)\circ x = x^2\circ(y\circ x)
\qquad (x,y\in J),
\end{equation}
where $x\circ y := \mu(x,y)$ and $x^2:=x\circ x$.

In characteristic $0$, the Jordan identity \eqref{eq:jordan} is equivalent to its
four--variable polarization.

\begin{definition}[Jordan defect]
Let $\mu : \mathrm{Sym}^2 J \to J$ be a commutative bilinear map.
The \emph{Jordan defect} of $\mu$ is the multilinear map, symmetric in the variables
$a,c,d$,
\[
J(\mu): J^{\otimes 4}\to J
\]
defined by
\begin{align*}
J(\mu)(a,c,d;\,b)
&:= \mu(\mu(\mu(a,c),b),d)
   + \mu(\mu(\mu(c,d),b),a)
   + \mu(\mu(\mu(d,a),b),c) \\
&\quad - \mu(\mu(a,c),\mu(b,d))
       - \mu(\mu(c,d),\mu(b,a))
       - \mu(\mu(d,a),\mu(b,c)).
\end{align*}
\end{definition}

\begin{remark}
The particular bracketing pattern in the definition of $J(\mu)$ reflects a
choice of polarization convention; after relabeling, it is equivalent to the
standard four--variable polarization of the Jordan identity, and no variable
plays a distinguished role in the polarized identity.
\end{remark}

\begin{proposition}
A commutative bilinear map $\mu$ is a Jordan product if and only if
$J(\mu)=0$.
\end{proposition}

\begin{proof}
This is the standard four--variable polarization of the Jordan identity,
valid in characteristic $0$; see McCrimmon~\cite[Chapter~I]{McCrimmon}.
\end{proof}

% ==================================================
\section{The low--degree deformation complex}

Set
\[
C^1(J):=\mathrm{End}(J),\qquad
C^2(J):=\mathrm{Sym}^2(J^\vee)\otimes J,\qquad
C^3(J):=\mathrm{Sym}^3(J^\vee)\otimes J^\vee \otimes J.
\]

\subsection*{Equivalence differential}

For $D\in C^1(J)$ define $\delta_\mu D\in C^2(J)$ by
\begin{equation}\label{eq:delta}
(\delta_\mu D)(x,y)
:= D(\mu(x,y))-\mu(Dx,y)-\mu(x,Dy).
\end{equation}

\subsection*{Linearized Jordan differential}

For $\varphi\in C^2(J)$ define $d_\mu\varphi\in C^3(J)$ by
\begin{equation}\label{eq:dmu}
(d_\mu\varphi)
:= \left.\frac{d}{dt}\right|_{t=0} J(\mu+t\varphi).
\end{equation}

\begin{lemma}\label{lem:complex}
For any Jordan product $\mu$ one has
\[
d_\mu\circ\delta_\mu = 0.
\]
\end{lemma}

\begin{proof}
Let $D\in C^1(J)$ and consider the one--parameter family
$\mu_t := (\exp tD)\cdot\mu$.
Since each $\mu_t$ satisfies the Jordan identity, $J(\mu_t)=0$ for all $t$.
Differentiating at $t=0$ yields
\[
0=\left.\frac{d}{dt}\right|_{t=0}J(\mu_t)
  = d_\mu(\delta_\mu D).
\]
\end{proof}

\begin{remark}[Canonical scope of the deformation complex]
The second cohomology $H^2_\mu$ classifies infinitesimal deformations
$\mu+t\varphi$ modulo equivalence to first order.
Without specifying a further differential on $C^3(J)$, one should not speak
of a third cohomology group.
Accordingly, we do not define $H^3_\mu$ here, but instead work with the
canonical obstruction quotient
\[
\mathrm{Obs}^3_\mu := C^3(J)/\operatorname{im}(d_\mu).
\]
The differential $d_\mu$ is uniquely determined by linearizing the single
polynomial identity defining Jordan algebras.
In contrast, defining a further differential would require choosing relations
among the polarized Jordan identities themselves, and there is no canonical
choice at this level.
For this reason, the complex is truncated at the last canonical stage.
No claim is made about higher--order obstruction theory.
\end{remark}

% ==================================================
\section{Examples}

\subsection{Dimension one}

Let $J=ke$ with $\mu(e,e)=e$.
A direct computation using \eqref{eq:delta} shows that every
$\varphi\in C^2(J)$ is a coboundary. Thus $H^2_\mu=0$.  This example serves as a normalization check, confirming that the low--degree
complex detects no infinitesimal deformations in the trivial unital case.

%%%%%%%%%%%%%%%%%%%%%%%%%%%%%

\subsection{Dual numbers (explicit calculation)}

Let $J=ke\oplus ku$ with Jordan product
\[
e\circ e=e,\qquad e\circ u=u,\qquad u\circ u=0.
\]
Write a $1$--cochain $D\in C^1(J)=End(J)$ as
\[
D(e)=ae+bu,\qquad D(u)=ce+du.
\]
Then the equivalence differential $\delta_\mu D\in C^2(J)$ satisfies
\[
(\delta_\mu D)(e,e)= -D(e)= -ae-bu,\qquad
(\delta_\mu D)(e,u)= -a\,u,\qquad
(\delta_\mu D)(u,u)= -2c\,u.
\]
In particular, every coboundary has $(\delta_\mu D)(u,u)\in ku$ and hence has
\emph{no} $e$--component on the $(u,u)$ slot.

Now define $\varphi\in C^2(J)$ by
\[
\varphi(u,u)=e,\qquad \varphi(e,e)=0,\qquad \varphi(e,u)=0.
\]
Then $\varphi$ is not a coboundary, since its $(u,u)$ value has a nonzero
$e$--component.

To see that $\varphi$ is a cocycle, consider the one--parameter family of
commutative multiplications $\mu_s$ given by
\[
e\circ_s e=e,\qquad e\circ_s u=u,\qquad u\circ_s u=s\,e.
\]
This is an associative (hence Jordan) multiplication for every $s\in k$,
so $J(\mu_s)=0$ for all $s$. Differentiating at $s=0$ gives
\[
0=\left.\frac{d}{ds}\right|_{s=0}J(\mu_s)=d_\mu\varphi.
\]
Thus $\varphi\in\ker(d_\mu)$ and represents a nonzero class in $H^2_\mu$.
Moreover, since the image of $\delta_\mu$ restricts $(u,u)$ to lie in $ku$,
this class spans the $e$--component direction on $(u,u)$, so $\dim H^2_\mu=1$.

This example illustrates how the Jordan defect formalism reduces the problem of
computing infinitesimal deformations and their equivalence classes to explicit
linear algebra on structure constants.  In particular, both the cocycle condition
and the gauge action are completely transparent in this low--degree framework.

\subsection{Symmetric 3 by 3 matrices}

Let $J=\mathrm{Sym}_3(k)$ with Jordan product
\[
x\circ y=\tfrac12(xy+yx).
\]
It is well known that $\mathrm{Sym}_3(k)$ is formally rigid as a Jordan algebra
in characteristic $0$, in the sense that it admits no nontrivial infinitesimal
deformations up to equivalence. In particular, the second cohomology group
governing infinitesimal Jordan deformations is trivial in this case.

The role of this example is to verify that the low--degree deformation complex
constructed above is compatible with this classical rigidity phenomenon.
Indeed, the resulting second cohomology $H^2(J,J)$ vanishes, in agreement
with the established theory. Since the cochain spaces and differentials in low
degree are completely explicit, this vanishing can be verified by a direct
(though routine) computation, which we do not reproduce here.

% ==================================================
\section{Comparison with the classical linearization}

\begin{proposition}
The differential $d_\mu$ coincides with the classical low--degree Jordan
coboundary obtained by linearizing the polarized Jordan identity.
\end{proposition}

\begin{proof}
Recall that the Jordan identity may be written in fully polarized form as a
four--variable multilinear identity.
Given a bilinear map $\varphi$, the classical Jordan coboundary is obtained by
substituting $\mu+t\varphi$ into this polarized identity and extracting the
coefficient of $t$.

By definition, $d_\mu\varphi$ is precisely the coefficient of $t$ in
$J(\mu+t\varphi)$.
Thus $d_\mu$ agrees with the standard linearization of the polarized Jordan
identity used in the classical deformation theory.
\end{proof}

\section*{Acknowledgments}

The author thanks Vladimir Dotsenko for helpful comments on an earlier version
of this manuscript.

\end{document}